\documentclass[journal]{IEEEtran}
\usepackage{amsmath,graphicx,amssymb,algorithm,algorithmic,multirow}
\newcommand{\tx}{\text}
\newtheorem{lemma}{Lemma}

\newtheorem{theorem}{Theorem}

\title{On-line Decentralized Charging of Plug-In Electric Vehicles in Power Systems}
\author{Qiao Li, \IEEEmembership{Student Member,~IEEE,} Tao Cui, \IEEEmembership{Student Member,~IEEE,} Rohit~Negi, \IEEEmembership{Member,~IEEE,}\\ Franz Franchetti, \IEEEmembership{Member,~IEEE,} and Marija D. Ili\'c, \IEEEmembership{Fellow, IEEE} \thanks{All authors are with the Department of Electrical and Computer Engineering, Carnegie Mellon University, Pittsburgh, PA, 15213 USA (email: \{qiaol, tcui, negi, franzf, milic\}@ece.cmu.edu).}}

\begin{document}
\maketitle

\begin{abstract}
  The concept of plug-in electric vehicles (PEV) are gaining increasing popularity in recent years, due to the growing societal awareness of reducing greenhouse gas (GHG) emissions, and gaining independence on foreign oil or petroleum. Large-scale deployment of PEVs currently faces many challenges. One particular concern is that the PEV charging can potentially cause significant impacts on the existing power distribution system, due to the increase in peak load. As such, this work tries to mitigate the impacts of PEV charging by proposing a decentralized smart PEV charging algorithm to minimize the distribution system load variance, so that a `flat' total load profile can be obtained. The charging algorithm is myopic, in that it controls the PEV charging processes in each time slot based entirely on the current power system states, without knowledge about future system dynamics. We provide theoretical guarantees on the asymptotic optimality of the proposed charging algorithm. Thus, compared to other forecast based smart charging approaches in the literature, the charging algorithm not only achieves optimality asymptotically in an on-line, and decentralized manner, but also is robust against various uncertainties in the power system, such as random PEV driving patterns and distributed generation (DG) with highly intermittent renewable energy sources.
\end{abstract}


\section*{Nomenclature}
\addcontentsline{toc}{section}{Nomenclature}
\begin{IEEEdescription}[\IEEEusemathlabelsep\IEEEsetlabelwidth{$V_1,V_2,V_3$}]
  \item[$A_i(n)$] Energy consumption of PEV $i$ at time slot $n$.
  \item[$A^{\max}_i$] Maximum energy consumption of PEV $i$ during a time slot.
  \item[$C_i$] Charging reference offset for PEV $i$.
  \item[$f(d)$] Daily average charging cost for day $d$.
  \item[$f^{\max}$] Maximum charging cost in a time slot.
  \item[$N$] Number of customers in the system.
  \item[$P_i(n)$] Charging power of PEV $i$ at time slot $n$.
  \item[$P^{\max}_i$] Maximum charging power of PEV $i$.
  \item[$S_{\tx{net}}(n)$] Total net base load observed by the substation at time slot $n$, which may include the output of distributed generators as \emph{negative} load.
  \item[$S_{\tx{net}}^{\tx{max}}$] Maximum total net base load in a time slot.
  \item[$\Delta t$] Length of a time slot.
  \item[$T$] Number of time slots for one day.
  \item[$U_i(n)$] Energy queue length of PEV $i$ at time slot $n$.
  \item[$U^{\tx{ref}}(n)$] PEV charging reference signal set by the aggregator at time slot $n$.
  \item[$\beta$] Weight of the charging cost.
  \item[$\eta_i$] Charging efficiency of PEV $i$.
  \item[$\mu(d)$] Daily average of the total net load for day $d$.
  \item[$\chi_i(n)$] Indicator function of whether PEV $i$ is available for charging at time slot $n$.  
\end{IEEEdescription}

\section{Introduction}
\label{sec_intro}

\IEEEPARstart{T}{he} growing societal awareness of environmental issues, as well as ongoing concerns about reducing the dependence on foreign oil or petroleum, have made the concept of plug-in electric vehicles (PEV) very popular during the past few years \cite{boulanger11}. PEVs, such as electric vehicles (EV) and plug-in hybrid electric vehicles (PHEV), can contribute to resolving these energy security issues by reducing the greenhouse gas (GHG) emissions and petroleum consumption in the transportation sector. Currently, large-scale implementation of PEVs in the near future is being planned. For example, the present US administration has planned 1 million PHEVs by 2015. Realizing the critical role of PEV in the auto industry, many automakers, such as Toyota, Nissan, GM, BYD, Fisker and Tesla are planning to produce dozens of types of PEVs \cite{boulanger11}, with many coming out starting by 2012.

It has been widely recognized that high penetration level of PEVs will cause significant impacts on the existing power system, in particular at the distribution level \cite{heydt83, meliopoulos09, clement-nyns10, fernandez11, lopes11, dyke10, schneider08, wu11}. Without proper coordination, it will be very likely that most of these PEVs will start charging during the overall peak load period \cite{heydt83}, causing severe branch congestions and voltage problems. Some studies \cite{lopes11, schneider08} have shown that the existing distribution system infrastructure may only support a very low PEV penetration level (such as $10\%$) without grid operation procedure changes or additional grid infrastructure investments. On the other hand, studies have also shown the promising result of mitigating the impact of PEV charging by \emph{coordinated charging}, which can effectively shifts the PEV load to the off-peak period. For example, it has been shown that coordinated charging can achieve around $50\%$ PEV penetration level in certain existing distribution systems \cite{schneider08}. Thus, it is crucial to design effective coordinated PEV charging algorithms for large-scale deployment of PEVs in the current power system, in order to achieve efficient grid operation, as well as bypassing or deferring the costly grid infrastructure investment. In other words, the existing `cyber' infrastructure may be used to resume the `physical' infrastructure in the power cyber-physical system (CPS).

This work contributes to the integration of the PEV into the power system by proposing an on-line smart charging algorithm. The algorithm tries to minimize the distribution system \emph{load variance}. The reason that such an objective function is chosen is as follows. Firstly, the minimum load variance objective function can achieve a perfect `valley-filling' charging profile, in the sense that at the optimal solution, the total load profile is as flat as it can possibly be \cite{gan11, ma10}. This implies that the PEV charging load can be efficiently `spread' among the off-peak periods, which may help achieve a higher PEV integration level in the existing power system, as well as lowering the distribution system loss \cite{sortomme11}, as compared to the other smart charging algorithms, in particular the ones based on electricity price \cite{lopes11, rotering11}. In the latter case, it is possible that a new `PEV charging peak' can form during the midnight, as many PEVs start charging simultaneously, as triggered by low electricity price. Secondly, the minimum load variance formulation is convex, and therefore can be solved exactly and quickly. Lastly, the proposed objective function is very flexible, which can also be used for other applications, such as unidirectional vehicle-to-grid (V2G) ancillary services. This can be done by `modulating' the base load according to external regulation signals. The details will be discussed later in this paper. 


As the main contribution of this work, this paper proposes an \emph{on-line decentralized} algorithm to solve the minimum load variance PEV charging problem. The charging algorithm is easy to implement. In each time slot, `binary' charging decisions (charge or not charge) are made \emph{locally} by each vehicle after comparing its battery's state of charge (SoC), which is the available percentage of the battery capacity, to a charging reference signal set by an aggregator, which is owned by the utility. The charging reference is carefully chosen, based entirely on the \emph{current states} of the power system, such as the SoC values of the plugged-in PEVs, the output of distributed generation (DG), the household base loads, as well as external regulation signals. In particular, the charging reference signal is chosen to greedily optimizes a function of the current system state, maximizing a battery `energy queue' weighted charging power, penalized by a quadratic function of total load in each time slot. Somewhat surprisingly, we will show that such a \emph{myopic} charging algorithm achieves \emph{the same asymptotic performance as compared to any optimal scheduling algorithm with perfect one-day ahead forecast of all uncertain parameters}. In other words, our algorithm achieves the same optimality as compared to conventional approaches, such as dynamic programming \cite{bertsekas95}, while \emph{dramatically reducing the computational complexity}. In fact, the proposed charging algorithm is closely related to the celebrated optimal \emph{max-weight} policy in the stochastic control literature \cite{meyn08, neely10}. Such max-weight type algorithms are intimately related to the stochastic sub-gradient algorithm \cite{neely10} and dynamic programming \cite{meyn08}, and has found numerous successes in diverse areas such as computer networks \cite{mckeown99}, wireless networks \cite{li11}, \cite{li09} and power transmission system \cite{meyn08}. In this paper, we extend the application of the max-weight algorithm to the PEV charging problem in the power distribution system, and prove the optimality results.

The smart PEV charging is an emerging area of research, which has been subject to continuing investigations. While there are many algorithms addressing smart PEV charging, most of them are solved in a centralized manner \cite{clement-nyns10, lopes11, sortomme11, saber11, sundstroem10}, assuming sufficient accuracy on day-ahead predictions about the stochastic dynamics in the power system. The real-time, decentralized charging issues are addressed by very few works, such as \cite{deilami11, turitsyn10}, where it is very challenging to obtain performance guarantees. Note that the forecast based algorithms are vulnerable to the prediction errors. In particular, they may face significant performance degradation as the current distribution system gradually evolves into the future \emph{smart distribution system} \cite{lasseter11}, where large-scale integration of distributed generation with intermittent renewable sources can result in highly uncertain stochastic dynamics, which are hard to predict from one day ahead. On the other hand, the proposed PEV charging algorithm in this paper does not suffer from such performance loss, as the PEV charging decisions are determined completely in real time, where power system states can be observed with very good accuracy. Thus, the proposed PEV charging algorithm can achieve both \emph{optimality} and \emph{robustness} in an on-line, decentralized manner.


%

The rest of this paper is organized as follows. In Section \ref{sec:model}, we introduce the model and formulate the optimal PEV charing problem. Section \ref{sec:charging} proposes and analyzes the optimal PEV charging algorithm, and Section \ref{sec:simulation} demonstrates the simulation results. Finally, Section \ref{sec:conclusion} concludes this paper.

\section{Assumptions and Modeling}
\label{sec:model}

In this section, we introduce the system model and formulate the minimum load variance PEV charging problem.

\subsection{Battery State Model}

A discrete-time system is considered in this paper, where the length of each time slot matches the typical sampling and operation time scale of the PEV aggregation unit. For example, the sampling rate of the distribution system load can be on a 15-min time basis \cite{clement-nyns10}. The PEV aggregator is assumed to be owned by a utility company, which is concerned about potential grid problems with PEV charging. Thus, the goal of the PEV aggregator is to achieve on-line, decentralized coordination of PEV charging, so as to minimize its impact on the distribution system.

Denote $T$ as the total number of time slots for each day. We assume that there are $N$ customers in the power system, and denote $ S_{\text{net}}(n)$ as the total \emph{net base load} as observed by the aggregator during time slot $n$, which may include the output power from distributed generation as \emph{negative loads}:
\begin{equation}
  S_{\tx{net}}(n) = S_{\tx{base}}(n)-{S}_{\tx{DG}}(n)
\end{equation}
where $S_{\tx{base}}$ and $S_{\tx{DG}}$ are the total base load and total distributed generation, respectively. Denote $U_i(n)$ as the `energy queue length' at PEV $i$ at the beginning of time slot $n$, which is the amount of energy that needs to be charged to refill the battery of PEV $i$. That is,
\begin{equation}
  U_i(n)= (1-{SoC}_i(n))\times {Cap}_i
\end{equation}
where ${SoC}_i$ and ${Cap}_i$ are the state of charge and capacity of PEV $i$'s battery, respectively. Thus, the queueing dynamics of the $U_i(n)$ can be expressed as the following:
\begin{equation}
  U_i(n) = U_i(n-1) - \eta_i P_i(n)\Delta t + A_i(n)
  \label{eqn:q_ev}
\end{equation}
In above, $P_i(n)$ is the charging power of PEV $i$ during time slot $n$, $\Delta t$ is the slot length, and $\eta_i$ is the battery charging efficiency of the PEV $i$, which depends on the PEV charger and battery types. It is assumed that
\begin{equation}
  0\leq P_i(n)\leq P^{\max}_i, 1\leq i\leq N \label{eqn:sigma_bound}
\end{equation}
where $P_i^{\max}$ depends on the charging circuit rating for PEV $i$. For example, for the standard 120-V/16-A wall outlet, the maximum charging power is 1.92kW. One technical issue is that the `energy queue length' cannot be negative, which implies that
\begin{equation}
  \eta_iP_i(n)\Delta t\leq U_i(n-1), 1\leq i\leq N \label{eqn:pi_pimax}
\end{equation}
This simply specifies that a battery will not further charge once it is full. $A_i(n)$ is the \emph{random} energy consumption of the PEV during time slot $n$. Bidirectional V2G applications are not considered in this paper, and will be addressed in future research. In this paper, the energy consumption $A_i(n)$ can be only caused by the \emph{random} PEV driving activity \cite{wu11} at time slot $n$, during which the vehicle is not able to draw power from the distribution system. Such charging availability is formally described by $\chi_i(n)$, which is an indicator function of whether PEV $i$ is plugged into the power grid at time slot $n$, i.e., $\chi_i(n)=1$ if PEV $i$ is plugged into the power grid at time slot $n$, otherwise $\chi_i(n)=0$. Thus, we have $P_i(n)=0$ if $\chi_i(n)=0$, since the PEV can not draw any power from the grid when it is not plugged into the power grid.

\subsection{Minimum Load Variance PEV Charging}

This paper focuses on coordinated PEV charging to minimize the total load variance as seen by the PEV aggregator. The optimal coordinated PEV charging problem can be formulated as follows:
\begin{IEEEeqnarray}{rCl}
  &\min_{\{P_i(n)\}}& f={1\over T}\sum_{n=1}^T (S_{\text{net}}(n)-\mu+\sum_{i=1}^N P_i(n))^2\label{eqn:cost}\\
  &\text{s.t. }& (\ref{eqn:q_ev}), (\ref{eqn:sigma_bound}), (\ref{eqn:pi_pimax}) \label{eqn:cost_aggreg}\\
  && P_i(n)=0 \text{ if }\chi_i(n)=0, \forall i, n\label{eqn:driving_pattern}\\
  && U_i(T)=U_i(0), 1\leq i\leq N\label{eqn:stability}
  \end{IEEEeqnarray}
In above, the cost function $f$ in (\ref{eqn:cost}) corresponds to the total load variance as seen by the PEV aggregator, where
\begin{equation}
  \mu = {1\over T}\sum_{n=1}^T(S_{\text{net}}(n)+\sum_{i=1}^N {A_i(n)\over \eta_i \Delta t})
  \label{eqn:def_mu}
\end{equation}
is the average total load during one day. The constraint in (\ref{eqn:stability}) essentially requires that all PEV energy consumptions during the day should be met by the PEV charging schedules $\{P_i(n)\}$. Note that our formulation has the same objective function as the minimum load variance formulation in \cite{sortomme11}, with more detailed modeling of the battery states and PEV driving processes as in (\ref{eqn:q_ev}) and (\ref{eqn:driving_pattern}). Further, note that both sequences of the net base load $\{S_{\tx{net}}(n)\}$ and the PEV driving variables $\{\chi_i(n)\}$, $\{A_i(n)\}$ are \emph{stochastic processes}. In the literature, such randomness is treated either as deterministic \cite{rotering11}, or by stochastic programming methods \cite{clement-nyns10}, which require the knowledge of the joint probability distribution of correlated random variables.

One obstacle in achieving completely myopic, real-time PEV charging is that the minimum load variance in (\ref{eqn:cost})-(\ref{eqn:stability}) requires knowledge of the average load $\mu$, which is a function of PEV loads during the whole day. On the other hand, such average load information is not needed to solve (\ref{eqn:cost})-(\ref{eqn:stability}), as we show in the following theorem:

\begin{theorem}
  Consider the following optimization problem:
  \begin{IEEEeqnarray}{rCl}
    &\min_{\{P_i(n)\}}& \tilde{f}={1\over T}\sum_{n=1}^T (S_{\text{net}}(n)+\sum_{i=1}^N P_i(n))^2\label{eqn:cost2}\\
  &\text{s.t. }& (\ref{eqn:cost_aggreg}), (\ref{eqn:driving_pattern}), (\ref{eqn:stability})\nonumber
  \end{IEEEeqnarray}
Any optimal solution of (\ref{eqn:cost2}) is also optimal for (\ref{eqn:cost})-(\ref{eqn:stability}).
  \label{theorem:equiv}
\end{theorem}

\begin{IEEEproof}
  See Appendix \ref{apdx:proof_equiv}.
\end{IEEEproof}
Thus, in order to solve the original minimum variance PEV charging problem, it is sufficient to solve (\ref{eqn:cost2}) instead, which does not include the average load $\mu$, and therefore allows easy implementation of myopic algorithms. Throughout of the rest of the paper, we will directly focus on solving (\ref{eqn:cost2}), bearing in mind that the original problem (\ref{eqn:cost})-(\ref{eqn:stability}) is also solved, according to Theorem \ref{theorem:equiv}.

In the literature, most smart charging algorithms \cite{clement-nyns10, lopes11, gan11, sortomme11} solve the optimization problem (\ref{eqn:cost2}) directly, assuming that all the random variables $\{S_{\tx{net}}(n)\}$, $\{A_i(n)\}$ and $\{\chi_i(n)\}$ (or their joint probability distribution) are known. Thus, the performance of these algorithms depends crucially on the prediction accuracy of these parameters, requiring knowledge of automotive driving patterns and base loads. While the former is clearly unrealistic, even loads may be unpredictable, if the system has a large number of small-scale distributed generators, which can use highly intermittent renewable energy sources, such as wind and solar. Prediction errors can cause substantial performance loss in these charging algorithms. As an alternative approach, in this paper, we propose a myopic PEV charging algorithm to solve (\ref{eqn:cost2}), which can achieve the optimal charging cost asymptotically, without suffering from day-ahead prediction errors.



%
\section{Decentralized Smart PEV Charging}
\label{sec:charging}

In this section, we describe and analyze the performance of the on-line decentralized PEV charging algorithm.

\subsection{Decentralized PEV Charging}

\begin{algorithm}
\caption{\bf On-line Decentralized PEV Charging}
\label{alg:pev}
\begin{algorithmic}[1]
  \STATE {\bf Initialization:} At the beginning of each time slot, the controller of each plugged-in PEV initializes its charging power $P_i\gets 0$, and the aggregator initializes the global charging reference as follows:
  \begin{equation}
    U_{\text{ref}} \gets (U_{\tx{ref}}^{\min}+U_{\tx{ref}}^{\max})/2
  \end{equation}
  where $U_{\tx{ref}}^{\min}$ and $U_{\tx{ref}}^{\max}$ are lower and upper bounds, respectively, which are selected based on historical data.
  \STATE Each vehicle controller updates its charging power:
    \begin{equation}
      P_i\gets\left\{
      \begin{array}{ll}
	0 & U_i(n-1)\leq {U_{\tx{ref}}/(\eta_i\Delta t)}-C_i\\
	P_i^{\max} & \tx{otherwise}
      \end{array}\right.
      \label{eqn:Pi}
    \end{equation}
    where $C_i$ is a properly chosen constant depending on the charging service contract of PEV $i$.

    \STATE Each vehicle controller submits the intended charging power $P_i$ to the aggregator, which sets $U_{\tx{ref}}^{\max}\gets U_{\tx{ref}}$ if
    \begin{equation}
      {U_{\tx{ref}}}> 2\beta(S_{\tx{net}}(n)+\sum_{i=1}^N {P_i})
      \label{eqn:nu}
    \end{equation}
    Otherwise, it sets $U_{\tx{ref}}^{\min}\gets U_{\tx{ref}}$. $\beta$ is a properly chosen constant to achieve a trade-off between refilling the batteries quickly and minimizing the charging cost.

    \STATE The charging reference $U_{\tx{ref}}$ is broadcast to each vehicle controller. If $|U_{\tx{ref}}^{\min}-U_{\tx{ref}}^{\max}|<\varepsilon'$, stop. Otherwise go to Step 2.
\end{algorithmic}
\end{algorithm}

\begin{figure}[!t]
   \begin{center}
     \includegraphics[width=3.5in]{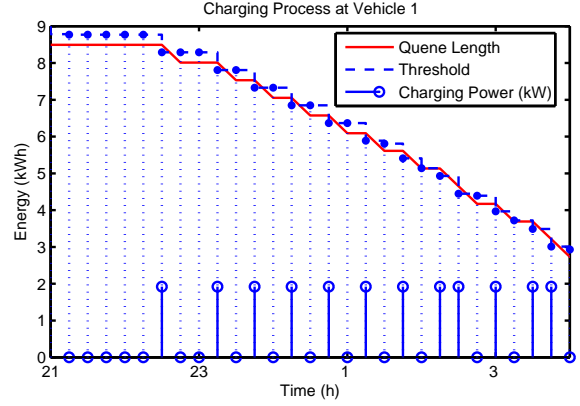}
     \caption{Illustration of the decentralized charging profile of a typical PEV.}
     \label{fig_charging_detail}
   \end{center}
 \end{figure}

The algorithm is shown in Algorithm \ref{alg:pev}. According to the algorithm, the charging decisions are made by each vehicle \emph{locally} by comparing its battery `energy queue length' (or equivalently, its SoC) to a charging reference set by the aggregator. Fig. \ref{fig_charging_detail} illustrates such local charging process at a PEV in the power system. In this way, the PEV charging processes are properly coordinated, so that the PEVs with higher `charging pressure' $U_i(n-1)$ (equivalently, lower SoC) can be served first. Further, compared to the centralized charging scheduling approaches, the algorithm is also less \emph{intrusive}, since the PEVs are the decision makers, and are only required to submit their charging power $\{P_i\}$ to the aggregator. Note that neither the charging decisions in (\ref{eqn:Pi}) nor the charging reference $U_{\tx{ref}}$ update is done in an ad hoc fashion. In fact, we will show the iterative procedures in Algorithm \ref{alg:pev} solve a generalized `max-weight' optimization problem.

\begin{theorem}
  In each time slot $n$, the charging power $\{P_i(n)\}$ computed by Algorithm \ref{alg:pev} solves the following optimization problem:
  \begin{IEEEeqnarray}{rCl}
  & \max_{\{P_i\}} & \sum_{i=1}^N (U_i(n-1)+C_i)\eta_i P_i\Delta t- \beta (S_{\text{net}}(n)+\sum_{i=1}^N P_i )^2\nonumber\\
  & \text{s.t. } & 0\leq P_i\leq P_i^{\max}, 1\leq i\leq N\nonumber\\
  && P_i\eta_i\Delta t\leq U_i(n-1), 1\leq i\leq N\nonumber\\
&& P_i = 0 \text{ if } \chi_i(n)=0, 1\leq i\leq N \label{eqn:max_weight}
\end{IEEEeqnarray}
For the case that the battery of PEV $i$ is fully charged in the middle of a time slot, $P_i(n)$ is interpreted as the `average charging power' $P_i(n)={U_i(n-1)/(\eta_i \Delta t)}$ for that particular time slot.
  \label{theorem:decentr}
\end{theorem}
\begin{IEEEproof}
  See Appendix \ref{apdx:proof_decentr}.
\end{IEEEproof}

We will continue with the above analysis in Section \ref{sec:guarantee} and present performance guarantees of Algorithm \ref{alg:pev}. Before that, we will discuss the intuition behind the algorithm, and its key properties.

\subsection{Discussions on the Charging Algorithm}

\subsubsection{Intuition behind the Algorithm}

Note that at the convergence, the charging reference can be written as $U_{\tx{ref}}=2\beta(S_{\tx{net}}(n)+\sum_{i=1}^N P_i)$. Thus, intuitively, during the peak period, the aggregator will set a high charging reference $U_{\text{ref}}$, so that the PEVs with low battery `energy queue length' (correspondingly, high SoC) will not charge. Similarly, for the off-peak period, the aggregator will set a low charging reference $U_{\text{ref}}$, so that the PEVs will `fill the valley'. Such an approach can be easily used to absorb generations with highly intermittent renewable sources, by treating them as negative loads in $S_{\tx{net}}(n)$. Note that it is also possible to `modulate' the net base load $S_{\tx{net}}(n)$ according to certain external regulation signals, in order to achieve V2G ancillary services. Detailed discussions on the specific methods, on the other hand, is out of the scope of this paper.

\subsubsection{Generalized Max-Weight Policy}

Anther way to interpret the algorithm is to inspect the optimization in (\ref{eqn:max_weight}). Note that the objective function in (\ref{eqn:max_weight}) has two parts: the `energy queue length weighted charging' term $(U_i(n-1)+C_i)\eta_i P_i\Delta t$, and the `charging cost' term $(S_{\text{net}}(n)+\sum_{i=1}^N P_i)^2$, which is weighted by parameter $\beta$. The first term is related to the `stability' requirement of the energy queues $\{U_i(n)\}$, which encourages the PEVs to charge at high power, in particular the ones with low SoC. On the other hand, the PEV charging is also penalized by the  second quadratic term, which is the instantaneous charging cost at time slot $n$, so as to prevent grid problems due to aggressive PEV charging. Finally, by adjusting the parameter $\beta$, one can achieve a trade-off between quickly charging the PEV batteries (maximizing the first term), and achieving minimum charging cost (minimizing the second term). The effect of $\beta$ will be illustrated in case studies. Note that $\beta$ needs to be chosen by the aggregator carefully, according to system specifications and historical data. 

\subsubsection{On-Off Charging} 

Algorithm \ref{alg:pev} has an interesting `on-off' charging property. That is, each PEV either charges at the maximum power, or do not charge at all. This can be seen from the algorithm specification in (\ref{eqn:Pi}), as well as in Fig. \ref{fig_charging_detail}. Thus, compared to the smart charging algorithms, which adjust charging powers continuously, the on-line decentralized PEV charging simplifies PEV charging circuit design. Similar type `on-off' charging profiles are also considered in \cite{han10}. 

\subsubsection{Communication Structure} 

\begin{figure}[t]
  \begin{center}
    \includegraphics[width=3.5in]{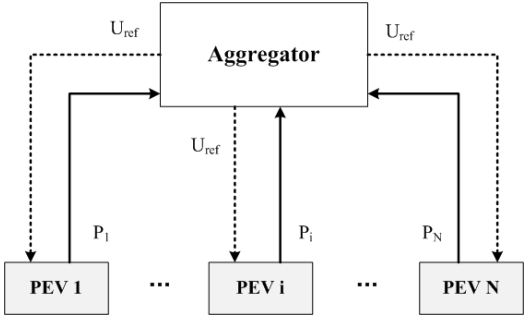}
  \end{center}
  \caption{The communication structure of the decentralized PEV charging algorithm.}
  \label{fig:comm}
\end{figure}

The communication structure for the algorithm is shown in Fig. \ref{fig:comm}. According to this structure, the vehicle controllers directly communicate with the aggregator, perhaps through extra intermediate levels of aggregation \cite{lopes11}. There is no direct communication between the PEVs. Further, the aggregator only needs a scalar, namely the aggregated PEV load $\big(\sum_{i=1}^N P_i\big)$, and broadcast back the charging reference $U_{\tx{ref}}$ to the PEVs. Therefore, the proposed scheme requires low communication capability, and is easy for real-time implementation.

\subsection{Performance Guarantees}
\label{sec:guarantee}

This section continues with the analysis of Theorem \ref{theorem:decentr} and demonstrates the asymptotic optimality of the proposed PEV charging algorithm, by comparing it against the optimal solution of (\ref{eqn:cost2}) over a long time period. As a minor technical assumption, we assume that the PEV energy consumptions $\{A_i(n)\}$ are \emph{strictly feasible}, in the sense that there is a charging schedule, which satisfies (\ref{eqn:sigma_bound}), (\ref{eqn:driving_pattern}), and that
\begin{equation}
  \sum_{n=1}^T P_i(n)\eta_i\Delta t\geq \sum_{n=1}^T A_i(n)+\varepsilon T, 1\leq i\leq N
  \label{eqn:strict_feasible}
\end{equation}
for a small constant $\varepsilon>0$. This simply means that the duration that a PEV is plugged into the distribution system is `more than just enough' to refill its daily energy consumption. We are interested in the performance of the algorithm over a period of $D$ days, and denote $\tilde{f}^\star(d)$ as the optimal cost of (\ref{eqn:cost2}) for the $d$-th day. We have the following theorem:


\begin{theorem}
  Algorithm \ref{alg:pev} achieves the following asymptotic average charging cost:
  \begin{IEEEeqnarray}{rCl}
    &&\limsup_{D\rightarrow\infty}{1\over DT} \sum_{n=1}^{DT} (\sum_{i=1}^N P_i(n)+S_{\tx{net}}(n))^2\leq\nonumber\\
    &&\qquad \qquad\limsup_{D\rightarrow\infty} {1\over D}\sum_{d=1}^D \tilde{f}^\star(d)+{B_1\over\beta}+{B_3(T+1)\over 2\beta}\label{eqn:cost_avg}
  \end{IEEEeqnarray}
  and the following average energy queue lengths:
  \begin{IEEEeqnarray}{rCl}
    &&\limsup_{D\rightarrow\infty}{1\over DT} \sum_{n=1}^{DT}\sum_{i=1}^N U_i(n)\leq {(B_2+B_3)(T+1)\over 2\varepsilon}\nonumber\\
    &&\qquad\qquad\qquad\qquad\qquad{B_1\over \varepsilon}+{\beta \tilde{f}^{\max}\over \varepsilon}-\sum_{i=1}^N C_i\label{eqn:queue_avg}
  \end{IEEEeqnarray}
  where the constants are defined as follows: $  \tilde{f}^{\max} = (S_{\tx{net}}^{\max}+\sum_{i=1}^N P_i^{\max})^2$, $B_1=\sum_{i=1}^N (C_i+K_i)K_i$, $B_2=\sum_{n=1}^N\varepsilon K_i$, $B_3=\sum_{n=1}^NK^2_i$, and the constant $K_i = \max(A_i^{\max}, \eta_iP_i^{\max}\Delta t)$.
  \label{theorem:scheduling}
\end{theorem}

\begin{IEEEproof}
  See Appendix \ref{apdx:scheduling}.
\end{IEEEproof}

From the above theorem, it is clear that, by properly adjusting the parameter $\beta$, the long term average charging cost can be made arbitrarily close to the optimal, with an optimality gap on an order of $O(1/\beta)$, as shown in (\ref{eqn:cost_avg}), while the average queue lengths grow as $O(\beta)$, as shown in (\ref{eqn:queue_avg}). The asymptotic bound provided by the above analysis is not tight in general. On the other hand, it is possible to provide tighter guarantees, by assuming specific probabilistic distributions on the stochastic processes, such as PEV driving patterns, base load, and renewable generation. Such a model dependent analysis is out of the scope of this paper, and will be addressed in future research. 

While this section showed the theoretical optimality of the decentralized PEV charging algorithm, in the next section, we demonstrate its performance in simulation to illustrate this point.



\section{Case Study}
\label{sec:simulation}

In this section, we test the performance of the proposed algorithm, and compare it against other prediction based algorithms. The simulations are done on IEEE 37-bus system and IEEE 123-bus system ~\cite{kersting01}, with MATLAB on a Xeon X3460 CPU with 8GB of RAM.

\subsection{Simulation Setup}

\subsubsection{Base Load}
 
\begin{figure}[!t]
  \begin{center}
     \includegraphics[width=3.5in]{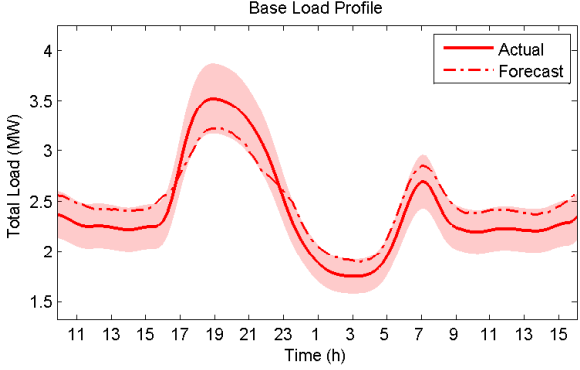}
     \caption{Base load profile used in the simulation with IEEE 37 bus system.}
     \label{fig:baseload37}
   \end{center}
 \end{figure}

The test systems are assumed to be residential distribution systems. Both systems use the same household base load curve, which is chosen according to the typical southern California residential load from the SCE website~\cite{SCEWebsite}, with proper scaling to match the national average household load of 1.3kW  \cite{EIALoad}. Fig. \ref{fig:baseload37} shows the daily curve of total base load used for the simulation with the IEEE 37-bus system. The IEEE 123-bus system uses a properly scaled curve with the same shape. For day-ahead prediction based PEV charging algorithms, we assume that they use forecast profiles such as the one in Fig. \ref{fig:baseload37}. Note that the shaded region corresponds to the 10\% mean average percentage error bounds, so that the errors of the forecast curve in Fig. \ref{fig:baseload37} are within $10\%$ of the actual load. For simplicity of comparison, we do not explicitly model DG in this simulation. But we want to emphasize that the prediction error will grow larger as the DG penetration level increases. Finally, the number of households connected to each bus will depend on the test feeder specification for each simulation scenario.

\subsubsection{PEV Specifications}

The PEV specifications and driving patterns are summarized in Table~\ref{tab:EVfacts}, which are obtained from typical PEV specifications~\cite{CheyVolt} and the national transportation survey~\cite{NHTS}. The daily consumption of 8.75kWh in Table~\ref{tab:EVfacts} is obtained from the 25 miles average daily commute distance and the PEV consumption rate of 34 kWh/100 miles. The average number of registered vehicle is assumed to be 1.8 per household \cite{NHTS}. In the simulation, The PEVs will be randomly assigned among all vehicles, according to the PEV penetration level. The random PEV driving pattern in this simulation is based on the national transportation survey \cite{NHTS}. A PEV is assumed to leave home at 7am, with zero-mean Guassian offset with standard deviation of 1 hour. Similarly, a PEV arrives home at 5pm, with zero-mean Gaussian offset with standard deviation of 2 hours. We focus on residential charging in this simulation, and assume that the PEV can only charge at home after it connects to the grid on returning home. 

\begin{table}[t]
\begin{center}
\caption{Vehicle Facts}
\begin{tabular}{cc}
\hline
Parameter & Value \\
\hline
Battery Capacity                 & 16 kWh   \\
Energy Usage per 100 miles             & 34 kWh   \\
Charging Rate (120 V, 16 A)      & 1.92 kW  \\
Average Daily Commute Distance   & 25 miles \\
Daily Consumption                & 8.75 kWh \\
Charging Efficiency & 0.90\\
\hline
\end{tabular}
\label{tab:EVfacts}
\end{center}
\end{table}

\subsection{IEEE 37-Bus System}

The PEV charging is first simulated in the standard IEEE 37-bus test feeder. In this case, the total number of vehicles is 3402. There are three types of smart charing algorithms considered in the simulation: 1) A static optimal charging algorithm which solves (\ref{eqn:cost2}), with \emph{perfect} knowledge of the day-ahead values of all randomness; 2) A static suboptimal charging algorithm, which solves (\ref{eqn:cost2}) using  \emph{imperfect forecast} of day-ahead load curve; and 3) On-line charging proposed by Algorithm \ref{alg:pev} in this paper. 

\begin{figure}[!t]
   \begin{center}
     \includegraphics[width=3.5in]{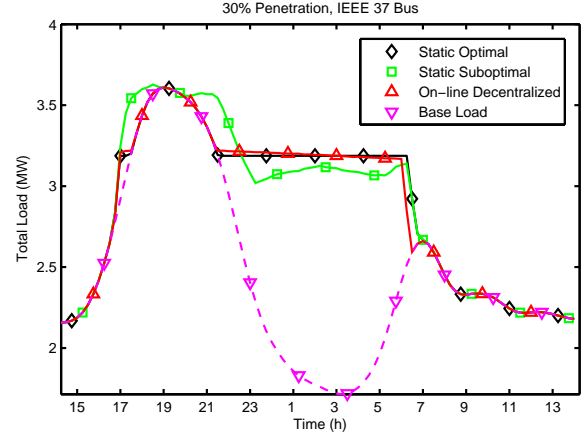}
     \caption{The total system loads with 30\% PEV penetration in the IEEE 37-bus system.}
     \label{fig:charging3730}
   \end{center}
 \end{figure}

\begin{figure}[!t]
   \begin{center}
     \includegraphics[width=3.5in]{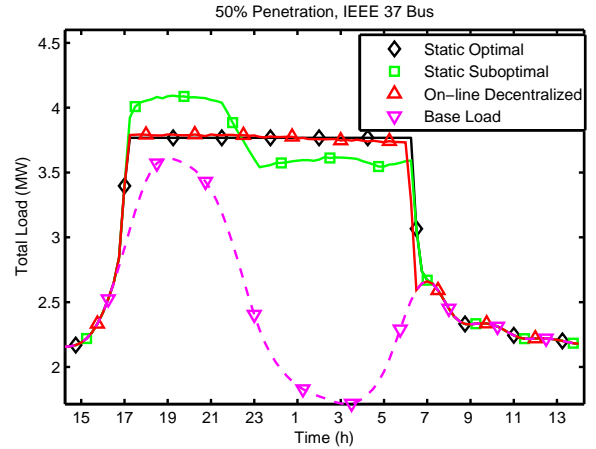}
     \caption{The total system loads with 50\% PEV penetration in the IEEE 37-bus system.}
     \label{fig:charging3750}
   \end{center}
 \end{figure}

The charging algorithms are simulated at PEV penetration levels of 30\% and 50\%. For the 30\% penetration case, $\beta = 0.0205$, and $C_i = 577$ for each vehicle, whereas for the 50\% penetration case, $\beta = 0.0161$, and $C_i = 534$ for each vehicle. In both cases, convergence can be observed after around 20 iterations for each time slot computation. The maximum total computation time of the on-line algorithm is 0.58 second for a 24-hour simulation scenario, while 3900 seconds for the static optimizations. Note the dramatic computation performance improvement for the case of on-line charging. This is due to the fact that each charging schedule is computed using \emph{only current system states}, which have much smaller dimension than the total state processes. In practice, the time scale of each time slot is on the order of minutes. Thus, the computation and communication requirement of the on-line charging algorithm can be easily satisfied. The results of total loads are shown in Fig. \ref{fig:charging3730} and Fig. \ref{fig:charging3750}, respectively. We have the following remarks.

\subsubsection{Valley Filling}

One can easily verify effectiveness of the minimum load variance formulation (\ref{eqn:cost}) by observing that, in both cases, the static optimal solution achieves a \emph{perfectly flat} total load curve at night. Thus, compared to other smart charging formulations, the ones based on electricity price in particular, the minimum load variance formulation can avoid the additional `midnight peak', which, in the extreme case, may cause similar grid congestion issues as uncoordinated charging. 

\subsubsection{Optimality}

The proposed on-line decentralized PEV charging achieves almost the same total load profile as the static optimal, \emph{even though the former does not need to know the driving pattern and loads in advance}. This further verifies the theoretical result in \ref{theorem:scheduling}. Thus, we can achieve the same performance as the static optimal, with much smaller computational overhead.

\subsubsection{Robustness}

The day-ahead prediction based algorithms are vulnerable to the forecast errors. This can be clearly observed from Fig. \ref{fig:charging3730} and Fig. \ref{fig:charging3750}, where the forecast based solutions can not achieve a flat profile in the presence of the load forecast error. In fact, we allowed these algorithms to know the exact driving patterns in advance, which is clearly unrealistic. On the other hand, the optimal decentralized charging algorithm is not affected by such forecast errors, since it is an on-line algorithm, which does not rely on forecasts.

\subsubsection{Service Differentiation} 

\begin{figure}[!t]
   \begin{center}
     \includegraphics[width=3.5in]{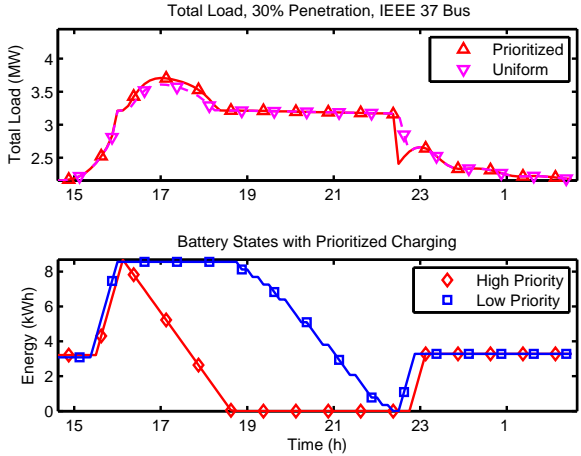}
     \caption{The total system load profile and typical battery state profile with 30\% PEV penetration in the IEEE 37-bus system with prioritized customers. Note that all customers in the `uniform' case have the same $C_i=577$.}
     \label{fig:battery3730}
   \end{center}
 \end{figure}

The constants $\{C_i\}$ are used to provide service differentiation among the PEVs, where a larger $C_i$ implies a higher priority for PEV $i$. That is, for two PEVs $i$ and $j$ with the same SoC, if $C_i>C_j$, PEV $i$ can start charging sooner than PEV $j$, according to (\ref{eqn:Pi}). In order to demonstrate this effect, we simulate the 30\% penetration case with two classes of customers, where $10\%$ customers have high priority, with $C_i=877$, and $90\%$ customers have low priority, with $C_i=577$. The total loads and battery states are shown in Fig. \ref{fig:battery3730}. One can clearly see that the energy queues of both classes are `stable'. Further, the high priority customers can finish charging much earlier than the customers with low priorities. Note that the values of $\{C_i\}$ have to be chosen carefully by the utility based on system specification and customers' choice, which is out of the scope of this paper.



  \begin{figure}[!t]
   \begin{center}
     \includegraphics[width=3.5in]{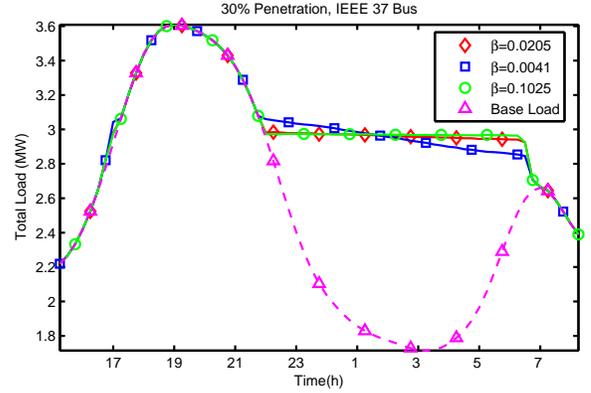}
     \caption{Sensitivity results of the total load with respect to changes in $\beta$ in the IEEE 37-bus system.}
     \label{fig:beta3730}
   \end{center}
 \end{figure}


\subsubsection{Effects of $\beta$} The parameter $\beta$ specifies the trade-off between charging the PEV batteries, and achieving a small charging cost. To illustrate this, we simulate the 30\% penetration case in Fig. \ref{fig:beta3730} with $\beta=0.0205, 0.0041,$ and $0.1025$, respectively. One can clearly observe that, when $\beta$ is small, the PEV charging is more greedy, as the weight of the charging cost is smaller. On the other hand, large $\beta$ will force the total load profile to be more flat, by penalizing the PEV charging power. In practice, $\beta$ has to be chosen carefully by the utility, using system specification and historical data.

 \subsection{IEEE 123-Bus System}

\begin{figure}[!t]
   \begin{center}
     \includegraphics[width=3.5in]{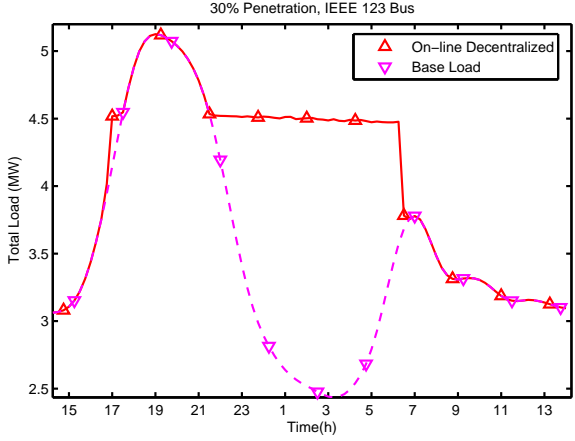}
     \caption{The total system load under the on-line decentralized algorithm with 30\% PEV penetration in the IEEE 123-bus system.}
     \label{fig:charging12330}
   \end{center}
 \end{figure}

\begin{figure}[!t]
   \begin{center}
     \includegraphics[width=3.5in]{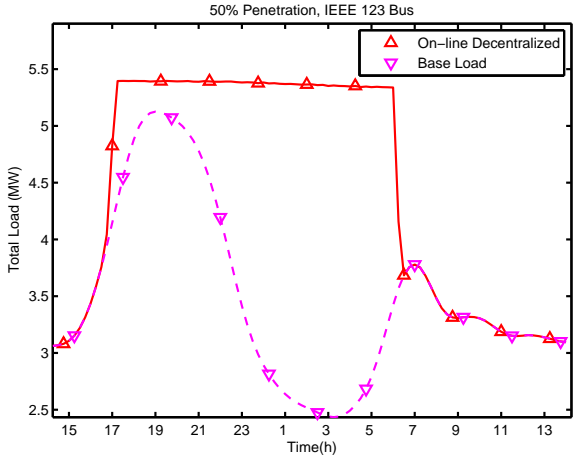}
     \caption{The total system load under the on-line decentralized algorithm with 50\% PEV penetration in the IEEE 123-bus system.}
     \label{fig:charging12350}
   \end{center}
 \end{figure}

We next simulate the PEV charging in IEEE 123-bus test feeder, where the total number of vehicles is 4832. In this case, one cannot implement the static optimizations on the machines specified at the beginning of this section, due to the large problem size (around 4600 thousand variables). On the other hand, \emph{the on-line charging can still be implemented}, as the charging decisions are made only using the current system states, which have much lower dimension (around 4.8 thousand variables). The charging algorithms are simulated at PEV penetration levels of 30\% and 50\%. For the 30\% case, $\beta=0.0205, C_i=815$, while for the 50\% case, $\beta=0.161,C_i=764$. In both simulations, convergence can be reached after around 20 iterations. The total computation time of the on-line algorithm is 0.75 second for a 24-hour simulation scenario. The results of total loads are shown in Fig. \ref{fig:charging12330} and Fig. \ref{fig:charging12350}, respectively. One can easily observe that, in both cases, the on-line algorithm can achieve a flat load profile. Thus, the on-line algorithm not only achieves the optimal charging cost, but also is dramatically computationally efficient, and much easier for implementation in large-scale systems, compared to other approaches forecast and static optimization based approaches.

\section{Conclusion}
\label{sec:conclusion}

This paper proposed an on-line decentralized PEV charging algorithm to minimize the distribution system total load variance. The charging decisions are made locally at each PEV based entirely on the power system states in each time slot. Since driving pattern and load forecasts are not needed in this algorithm, the charging performance is robust against various uncertainties in the system, as compared to the prediction based static optimization approaches in the literature. In the analysis, asymptotic results about the cost function and battery level stability are shown. The performance of the proposed charging algorithm is further compared against static smart charging algorithms by simulation results.

\appendices

\section{Proof of Theorem \ref{theorem:equiv}}
\label{apdx:proof_equiv}
\begin{IEEEproof}
  It is obvious that any feasible solution of (\ref{eqn:cost2}) is also feasible for (\ref{eqn:cost})-(\ref{eqn:stability}). Now, it is sufficient to show that the two objective functions are different by at most a constant that is independent of choice of charging schedules $\{P_i(n)\}$. Note that the constraint (\ref{eqn:stability}) implies that
  \begin{equation}
    \sum_{n=1}^T \eta_iP_i(n)\Delta t = \sum_{n=1}^T A_i(n)
  \end{equation}
  This, together with the definition of $\mu$ in (\ref{eqn:def_mu}) imply that
  \begin{equation}
    \mu = {1\over T}\sum_{n=1}^T (S_{\tx{net}}(n)+\sum_{i=1}^N P_i(n))
  \end{equation}
  Therefore, we have
    \begin{IEEEeqnarray}{rCl}
      f&=&\tilde{f}-{2\mu\over T}\sum_{n=1}^T (S_{\tx{net}}(n)+\sum_{i=1}^N P_i(n))+\mu^2\\
      &=&\tilde{f}-\mu^2,
    \end{IEEEeqnarray}
    where, by (\ref{eqn:def_mu}), $\mu$ is independent of the choice of schedules. Thus, the theorem holds.
\end{IEEEproof}

\section{Proof of Theorem \ref{theorem:decentr}}
\label{apdx:proof_decentr}
\begin{IEEEproof}
  We can formulate the optimization in (\ref{eqn:max_weight}) as the following equivalent problem:
  \begin{IEEEeqnarray}{rCl}
    & \max_{\{P_i\}, y} & \sum_{i=1}^N (U_i(n-1)+C_i) \eta_i P_i \Delta t - \beta y^2\label{eqn:mw_relaxed}\\
    & \text{subject to } & y = \sum_{i=1}^N P_i + S_{\text{net}}(n)\\
    && 0\leq P_i\leq P_i^{\max}, 1\leq i\leq N\label{eqn:pi_lim}\\
    && P_i\eta_i\Delta t\leq U_i(n-1), 1\leq i\leq N\label{eqn:pi_lim2}\\
    && P_i = 0 \text{ if } \chi_i(n)=0, 1\leq i\leq N\label{eqn:chi_lim}
  \end{IEEEeqnarray}
  Thus, one can easily write the dual problem as
   \begin{IEEEeqnarray}{rCl}
     & \min_{U_{\tx{ref}}}\max_{\{P_i, y\}} & \sum_{i=1}^N (U_i(n-1)+C_i) \eta_i P_i \Delta t - \beta y^2\nonumber\\
     &&\quad\qquad\qquad\qquad +U_{\tx{ref}}(y - \sum_{i=1}^N P_i - S_{\text{net}}(n))\nonumber\\
     & \text{subject to } & (\ref{eqn:pi_lim}), (\ref{eqn:pi_lim2}) \text{ and } (\ref{eqn:chi_lim})\label{eqn:dual1}
   \end{IEEEeqnarray}
   It can be easily shown that $y^\star = {U_{\tx{ref}}/ 2\beta}$. Thus, we can simplify (\ref{eqn:dual1}) as
  \begin{IEEEeqnarray}{rCl}
    & \min_{U_{\tx{ref}}}\max_{\{P_i\}} & \sum_{i=1}^N (U_i(n-1)+C_i-{U_{\tx{ref}}\over \eta_i\Delta t}) \eta_i P_i \Delta t\nonumber\\
    &&\qquad\qquad\qquad \qquad \qquad +{({U_{\tx{ref}}})^2\over 4 \beta}-U_{\tx{ref}}S_{\text{net}}(n)\nonumber\\
    & \text{subject to } & (\ref{eqn:pi_lim}), (\ref{eqn:pi_lim2}) \text{ and } (\ref{eqn:chi_lim})\label{eqn:dual2}
  \end{IEEEeqnarray}  
  One can immediately see that the decentralized charging in (\ref{eqn:Pi}) corresponds to the inner optimization problem with charging power constraint (\ref{eqn:pi_lim}), and the update on $U_{\tx{ref}}$ in Step 3 of Algorithm \ref{alg:pev} corresponds to standard binary search algorithm over the scalar $U_{\tx{ref}}$ to solve the convex dual problem.
\end{IEEEproof}

\section{Proof of Theorem \ref{theorem:scheduling}}
\label{apdx:scheduling}

In this section we prove Theorem \ref{theorem:scheduling}. Define a `Lyapunov' type function $F(n)$ as follows:
\begin{eqnarray*}
    F(n)={1\over 2}\sum_{i=1}^N (U_i(n)+C_i)^2+\beta\sum_{k=1}^n (S_{\tx{net}}(k)+\sum_{i=1}^N P_i(k))^2
\end{eqnarray*}
The key in proving both cost optimality as well as stability results involves analyzing the drift behavior of $F(n)$, which we will illustrate in the following. Note that similar techniques have also been applied in \cite{neely10} in the context of wireless networks.

\subsection{Cost Optimality}

The proof of (\ref{eqn:cost_avg}) requires several technical lemmas. The following lemma provided a bound on the single slot drift of $F(n)$.

\begin{lemma}
  The one-slot drift of $F(n)$ can be bounded as follows:
  \begin{IEEEeqnarray}{rCl}
  &&\Delta_1 F(0)\nonumber\\
  &=& F(1)-F(0)\nonumber\\
  &\leq &\sum_{i=1}^N (U_i(0)+C_i)(A_i(1)-\eta_i P_i(1)\Delta t)+B_3\nonumber\\
  &&\quad\qquad\qquad\qquad\qquad+\beta (S_\tx{net}(1)+\sum_{i=1}^N P_i(1))^2
\end{IEEEeqnarray}
\label{lem:single_drift}
\end{lemma}

\begin{IEEEproof}
For each PEV $i$, direct calculation shows that
\begin{eqnarray}
  &&{1\over 2}(U_i(1)+C_i)^2\nonumber\\
  &=&{1\over 2}\big(U_i(0)+A_i(1)-\eta_i P_i(1)\Delta t+C_i\big )^2\nonumber\\
&=&{1\over 2}(U_i(0)+C_i)^2+(U_i(0)+C_i)(A_i(1)-\eta_i P_i(1)\Delta t)\nonumber\\
&&\quad\qquad\qquad\qquad\qquad +{1\over 2}(A_i(1)-\eta_i P_i(1)\Delta t)^2
\end{eqnarray}
Thus, the following bound holds:
\begin{eqnarray}
  &&{1\over 2}(U_i(1)+C_i)^2-{1\over 2}(U_i(0)+C_i)^2\nonumber\\
  &=&(U_i(0)+C_i)(A_i(1)-\eta_i P_i(1)\Delta t)\nonumber\\
  &&\qquad\qquad\qquad\qquad\qquad+{1\over 2}(A_i(1)-\eta_i P_i(1)\Delta t)^2\nonumber\\
  &\leq &(U_i(0)+C_i)(A_i(1)-\eta_i P_i(1)\Delta t)+{1\over 2}K_i^2
\end{eqnarray}
where $K_i=\max(A^{\max}_i, \eta_i P_i^{\max}\Delta t)$. The last inequality holds because we assume that a PEV can either charge (in garage) or discharge (on the road) in each time slot, but not both. Thus, the lemma follows after summing up the above inequality over $N$ customers and adding the charging cost drift term $\beta (S_\tx{net}(1)+\sum_{i=1}^N P_i(1))^2$. 
\end{IEEEproof}

We next generalize the above bound to $T$ time slots.

\begin{lemma}
  Under Algorithm \ref{alg:pev}, the first $T$-slot drift of $F(n)$ can be bounded as
    \begin{IEEEeqnarray}{rCl}
      \Delta_T F(0) \leq B_1T+B_3{T(T+1)\over 2}+\beta \tilde{f}^{\star}(1)T
  \end{IEEEeqnarray}
  \label{lem:drift_Cost}
\end{lemma}

\begin{IEEEproof}
  Carrying out the procedure in Lemma \ref{lem:single_drift} over $T$ time slots, we obtain the following
  \begin{eqnarray}
    \Delta_T F(0)&\leq&\sum_{n=1}^T\sum_{i=1}^N (U_i(n-1)+C_i)(A_i(n)-\eta_i P_i(n)\Delta t)\nonumber\\
    &&\quad +B_3 T +\beta \sum_{n=1}^T(S_{\tx{net}}(n)+\sum_{i=1}^N P_i(n))^2\label{eqn:boundFnT1}
  \end{eqnarray}
  Now, we want to compare the bound on the RHS of (\ref{eqn:boundFnT1}) under the myopic charging schedule $\{P_i(n)\}$, which are computed by Algorithm \ref{alg:pev}, against a static optimal solution of (\ref{eqn:cost2}) with perfect knowledge of all stochastic dynamics, which we denote as $\{P^{\tx{stat}}_i(n)\}$. Because Algorithm \ref{alg:pev} always greedily maximizes (\ref{eqn:max_weight}), the RHS of the above can be further bounded as follows:
  \begin{eqnarray}
      &&\Delta_T F(0)\nonumber\\
      &\leq&\sum_{n=1}^T\sum_{i=1}^N (U_i(n-1)+C_i)(A_i(n)-\eta_i P_i^{\tx{stat}}(n)\Delta t)\nonumber\\
   &&(B_1+B_3) T +\beta \sum_{n=1}^T(S_{\tx{net}}(n)+\sum_{i=1}^N P_i^{\tx{stat}}(n))^2
   \label{eqn:driftT1}
  \end{eqnarray}
  The extra term $B_1=\sum_{i=1}^N (C_i+K_i)K_i$ is introduced to bound the case where the energy queues are small, so that the charging process stops in the middle of a time slot. In such a case, we have
  \begin{equation}
    (U_i(n-1)+C_i)\eta_iP_i^{\tx{stat}}(n)\Delta t\leq (K_i+C_i)K_i
  \end{equation}
  We next bound the drift in the energy queue lengths. Note that the energy queue length for each PEV $i$ can be written as the following:
\begin{equation}
  U_i(n)=U_i(0)+\sum_{k=1}^n \big(A_i(n)-\eta_iP_i(n)\Delta t\big)
\end{equation}
Thus, $U_i(n)$ can be bounded as follows:
\begin{equation}
  U_i(0)-K_i n \leq  U_i(n) \leq U_i(0)+K_i n\label{eqn:bound_Ui}
\end{equation}
which yields the following bound:
\begin{eqnarray}
&&U_i(n-1) \big(A_i(n)-\eta_iP_i^{\tx{stat}}(n)\Delta t\big) \leq\nonumber\\ 
&&\quad U_i(0)\big(A_i(n)- \eta_iP_i^{\tx{stat}}(n)\Delta t\big)+K^2_i (n-1) 
\end{eqnarray}
Plugging the above inequality into the (\ref{eqn:driftT1}) yields
\begin{eqnarray}
  &&\Delta_T F(0)\nonumber\\
  &\leq& \sum_{i=1}^N (U_i(0)+C_i)\Big(\sum_{n=1}^T\big(A_i(n)-\eta_i P_i^{\tx{stat}}(n)\Delta t\big) \Big)+ B_1 T \nonumber\\
  &&+B_3T+ B_3 \sum_{n=1}^T (n-1)+\beta \sum_{n=1}^T(S_{\tx{net}}(n)+\sum_{i=1}^N P_i^{\tx{stat}}(n))^2\nonumber\\
  &\stackrel{(a)}{=}& B_1 T+B_3{T(T+1)\over 2}+\beta \tilde{f}^{\star}(1)T
\end{eqnarray}
where step $(a)$ is because $\{P_i^{\tx{stat}}(n)\}$ solves the optimization (\ref{eqn:cost2}), and satisfies
\begin{eqnarray}
  \sum_{n=1}^T \big(A_i(n)-\eta_i P_i^{\tx{stat}}(n)\Delta t\big)&=&0\\
  {1\over T}\sum_{n=1}^T(S_{\tx{net}}(n)+\sum_{i=1}^N P_i^{\tx{stat}}(n))^2&=&\tilde{f}^\star(1)
\end{eqnarray}
Thus, the lemma follows.
\end{IEEEproof}

We next extend the above analysis from one day to $D$ days:
\begin{lemma}
  The drift of $F(n)$ over the first $DT$ time slots satisfies
    \begin{IEEEeqnarray*}{rCl}
      && \Delta_{DT} F(0)\leq B_1 DT + B_3{DT(T+1)\over 2}+\beta T\sum_{d=1}^D\tilde{f}^{\star}(d)
  \end{IEEEeqnarray*}
  \label{lem:drift_Cost_D}
\end{lemma}
\begin{IEEEproof}
  This can be simply obtained by summing the bound in Lemma \ref{lem:drift_Cost} over $D$ days.
\end{IEEEproof}

We are now ready to prove the first part of Theorem \ref{theorem:scheduling}.

\begin{IEEEproof}
  \emph{(Cost Optimality of Algorithm \ref{alg:pev})} According to the bound in Lemma \ref{lem:drift_Cost_D}, the average cost of the first $D$ days under Algorithm \ref{alg:pev} can be bounded as
  \begin{IEEEeqnarray}{rCl}
    &&{1\over DT} \sum_{n=1}^{DT}(S_{\tx{net}}(n)+\sum_{i=1}^N P_i(n))^2\nonumber \\
    &\leq& {1\over \beta DT} F(DT)\nonumber\\
    &=&{F(0)+\Delta_{DT}F(0)\over \beta DT }\nonumber\\
    &\leq& {F(0)\over DT } +{B_1\over \beta} + B_3 {T+1\over 2\beta}+ {1\over D}\sum_{d=1}^Df^\star(d)
  \end{IEEEeqnarray}
Thus, after taking $D\rightarrow\infty$, the cost optimality claim in (\ref{eqn:cost_avg}) follows. 
\end{IEEEproof}

\subsection{Stability}

The proof of the second part of Theorem \ref{theorem:scheduling} requires another set of lemmas bounding the drift of $F(n)$. We start with the first one the $T$-slot drift of the function $F(n)$.

\begin{lemma}
 Under Algorithm \ref{alg:pev}, the first $T$-slot drift of $F(n)$ can be bounded as  
  \begin{IEEEeqnarray}{rCl}
&&\Delta_T F(0) \leq -\varepsilon \sum_{n=1}^T\sum_{i=1}^N (U_i(n)+C_i) +B_1 T\nonumber\\
&&\qquad\qquad\qquad + (B_2+B_3){T(T+1)\over 2}+\beta \tilde{f}^{\max}T
  \end{IEEEeqnarray}
  \label{lem:boundFn2}
\end{lemma}
\begin{IEEEproof}
  Denote $\{P_i^{\tx{feas}}(n)\}$ as a sequence of \emph{strictly} feasible charging schedules, which satisfy (\ref{eqn:strict_feasible}). Because the charging schedules $\{P_i(n)\}$ by Algorithm \ref{alg:pev} always greedily maximize (\ref{eqn:max_weight}), we can obtain a similar bound as (\ref{eqn:driftT1}) as follows:
    \begin{eqnarray}
      &&\Delta_T F(0)\nonumber\\
      &\leq&\sum_{n=1}^T\sum_{i=1}^N (U_i(n-1)+C_i)(A_i(n)-\eta_i P_i^{\tx{feas}}(n)\Delta t)\nonumber\\
   &&\qquad +B_1T+B_3 T +\beta \sum_{n=1}^T(S_{\tx{net}}(n)+\sum_{i=1}^N P_i^{\tx{feas}}(n))^2\nonumber\\
   &\leq& \sum_{n=1}^T\sum_{i=1}^N (U_i(n-1)+C_i)(A_i(n)-\eta_i P_i^{\tx{feas}}(n)\Delta t)\nonumber\\
   &&\qquad\qquad\qquad\qquad+B_1T+B_3 T +\beta \tilde{f}^{\max}T\label{eqn:driftT2}
  \end{eqnarray}
  From (\ref{eqn:bound_Ui}), the following bound holds:
  \begin{eqnarray}
&&U_i(n-1) \big(A_i(n)-\eta_iP_i^{\tx{feas}}(n)\Delta t\big) \leq\nonumber\\ 
&&\quad U_i(0)\big(A_i(n)- \eta_iP_i^{\tx{feas}}(n)\Delta t\big)+K^2_i (n-1) 
\end{eqnarray}
Plugging the above into (\ref{eqn:driftT2}) and applying (\ref{eqn:strict_feasible}) yields
   \begin{eqnarray}
      &&\Delta_T F(0)\nonumber\\
      &\leq&\sum_{i=1}^N(U_i(0)+C_i)\Big(\sum_{n=1}^T (A_i(n)-\eta_i P_i^{\tx{feas}}(n)\Delta t)\Big)\nonumber\\
      &&\qquad\qquad\qquad\qquad +B_1T+B_3{T(T+1)\over 2}+\beta \tilde{f}^{\max}T\nonumber\\
      &\leq&-\varepsilon T\sum_{i=1}^N(U_i(0)+C_i)+B_3{T(T+1)\over 2}\nonumber\\
      &&\qquad\qquad\qquad\qquad\qquad\qquad+B_1T+\beta \tilde{f}^{\max}T
  \end{eqnarray}
  Finally, from (\ref{eqn:bound_Ui}) the above can be further bounded as
  \begin{eqnarray}
      &&\Delta_T F(0)\nonumber\\
      &\leq&-\varepsilon \sum_{i=1}^N\sum_{n=1}^T(U_i(n)- K_in+C_i)\nonumber\\
      &&\qquad\qquad\qquad+B_1T+B_3{T(T+1)\over 2}+\beta \tilde{f}^{\max}T\nonumber\\
      &=&-\varepsilon \sum_{i=1}^N\sum_{n=1}^T(U_i(n)+C_i)+B_2{T(T+1)\over 2}\nonumber\\
      &&\quad\qquad\qquad+B_1T+B_3{T(T+1)\over 2}+\beta \tilde{f}^{\max}T
  \end{eqnarray}
  from which the lemma follows.
\end{IEEEproof}

Similarly, the drift of $F(n)$ over $D$ days can be bounded by the following lemma.
\begin{lemma}
  The drift of $F(n)$ over the first $D$ days can be bounded as follows:
  \begin{IEEEeqnarray}{rCl}
    && \Delta_{DT} F(0)\leq -\varepsilon \sum_{n=1}^{DT}\sum_{i=1}^N (U_i(n)+C_i)+B_1DT\nonumber\\
    &&\qquad\qquad+ (B_2+B_3){DT(T+1)\over 2}+\beta D\tilde{f}^{\max}T
    \label{eqn:driftFnD2}
  \end{IEEEeqnarray}
  \label{lem:driftFnD2}
\end{lemma}
\begin{IEEEproof}
  This can be easily verified by applying Lemma \ref{lem:boundFn2} $D$ times and summing the bounds.
\end{IEEEproof}

We are now ready to prove the second part of Theorem \ref{theorem:scheduling}.
\begin{IEEEproof}
  \emph{(Stability of Algorithm \ref{alg:pev})} From (\ref{eqn:driftFnD2}) one can easily see that
  \begin{IEEEeqnarray}{rCl}
  &&{1\over DT} \sum_{n=1}^{DT} \sum_{i=1}^N (U_i(n)+C_i)\nonumber\\
  &\leq& {F(0)-F(DT)\over \varepsilon DT}+{B_1\over\varepsilon}+{(B_2+B_3)(T+1)\over 2\varepsilon}+ {\beta \tilde{f}^{\max}\over \varepsilon}\nonumber\\
&\leq& {F(0)\over \varepsilon DT}+{B_1\over\varepsilon}+{(B_2+B_3)(T+1)\over 2\varepsilon}+ {\beta \tilde{f}^{\max}\over \varepsilon}
\end{IEEEeqnarray}  
from which (\ref{eqn:queue_avg}) follows from above by taking $D\rightarrow\infty$.
\end{IEEEproof}

\bibliographystyle{IEEEtran}
\bibliography{IEEEabrv,ev}

\begin{thebibliography}{10}
\providecommand{\url}[1]{#1}
\csname url@samestyle\endcsname
\providecommand{\newblock}{\relax}
\providecommand{\bibinfo}[2]{#2}
\providecommand{\BIBentrySTDinterwordspacing}{\spaceskip=0pt\relax}
\providecommand{\BIBentryALTinterwordstretchfactor}{4}
\providecommand{\BIBentryALTinterwordspacing}{\spaceskip=\fontdimen2\font plus
\BIBentryALTinterwordstretchfactor\fontdimen3\font minus
  \fontdimen4\font\relax}
\providecommand{\BIBforeignlanguage}[2]{{%
\expandafter\ifx\csname l@#1\endcsname\relax
\typeout{** WARNING: IEEEtran.bst: No hyphenation pattern has been}%
\typeout{** loaded for the language `#1'. Using the pattern for}%
\typeout{** the default language instead.}%
\else
\language=\csname l@#1\endcsname
\fi
#2}}
\providecommand{\BIBdecl}{\relax}
\BIBdecl

\bibitem{boulanger11}
A.~Boulanger, A.~Chu, S.~Maxx, and D.~Waltz, ``Vehicle electrification: Status
  and issues,'' \emph{Proc. {IEEE}}, vol.~99, no.~6, pp. 1116--1138, Jun. 2011.

\bibitem{heydt83}
G.~Heydt, ``The impact of electric vehicle deployment on load management
  strategies,'' \emph{{IEEE} Trans. Power App. Syst.}, vol. PAS-102, no.~5, pp.
  1253--1259, May 1983.

\bibitem{meliopoulos09}
S.~Meliopoulos, J.~Meisel, G.~Cokkinides, and T.~Overbye, ``Power system level
  impacts of plug-in hybrid vehicles,'' Power Systems Engineering Research
  Center (PSERC), Tech. Rep., Oct. 2009.

\bibitem{clement-nyns10}
K.~Clement-Nyns, E.~Haesen, and J.~Driesen, ``The impact of charging plug-in
  hybrid electric vehicles on a residential distribution grid,'' \emph{{IEEE}
  Trans. Power Syst.}, vol.~25, no.~1, pp. 371--380, Feb. 2010.

\bibitem{fernandez11}
L.~Fernandez, T.~Roman, R.~Cossent, C.~Domingo, and P.~Frias, ``Assessment of
  the impact of plug-in electric vehicles on distribution networks,''
  \emph{{IEEE} Trans. Power Syst.}, vol.~26, no.~1, pp. 206--213, Feb. 2011.

\bibitem{lopes11}
J.~Lopes, F.~Soares, and P.~Almeida, ``Integration of electric vehicles in the
  electric power system,'' \emph{Proc. {IEEE}}, vol.~99, no.~1, pp. 168--183,
  Jan. 2011.

\bibitem{dyke10}
K.~Dyke, N.~Schofield, and M.~Barnes, ``The impact of transport electrification
  on electrical networks,'' \emph{{IEEE} Trans. Ind. Electron.}, vol.~57,
  no.~12, pp. 3917--3926, Dec. 2010.

\bibitem{schneider08}
K.~Schneider, C.~Gerkensmeyer, M.~Kintner-Meyer, and R.~Fletcher, ``Impact
  assessment of plug-in hybrid vehicles on pacific northwest distribution
  systems,'' in \emph{Proc. Power Energy Soc. Gen. Meet.}, Jul. 2008, pp. 1--6.

\bibitem{wu11}
D.~Wu, D.~Aliprantis, and K.~Gkritza, ``Electric energy and power consumption
  by light-duty plug-in electric vehicles,'' \emph{{IEEE} Trans. Power Syst.},
  vol.~26, no.~2, pp. 738--746, May 2011.

\bibitem{gan11}
L.~Gan, UT, and S.~Low, ``Optimal decentralized protocols for electric vehicle
  charging,'' in \emph{Proc. IEEE Conf. on Dec. and Contr.}, 2011.

\bibitem{ma10}
Z.~Ma, D.~Callaway, and I.~Hiskens, ``Decentralized charging control for large
  populations of plug-in electric vehicles: Application of the nash certainty
  equivalence principle,'' in \emph{Proc. IEEE Conf. on Dec. and Contr.}, Sep.
  2010, pp. 191--195.

\bibitem{sortomme11}
E.~Sortomme, M.~Hindi, S.~MacPherson, and S.~Venkata, ``Coordinated charging of
  plug-in hybrid electric vehicles to minimize distribution system losses,''
  \emph{{IEEE} Trans. Smart Grid}, vol.~2, no.~1, pp. 198--205, Mar. 2011.

\bibitem{rotering11}
N.~Rotering and M.~Ilic, ``Optimal charge control of plug-in hybrid electric
  vehicles in deregulated electricity markets,'' \emph{{IEEE} Trans. Power
  Syst.}, vol.~26, no.~3, pp. 1021--1029, Aug. 2011.

\bibitem{bertsekas95}
D.~Bertsekas, \emph{Dynamic programming and optimal control}.\hskip 1em plus
  0.5em minus 0.4em\relax Athena Scientific, 1995.

\bibitem{meyn08}
S.~Meyn, \emph{Control techniques for complex networks}.\hskip 1em plus 0.5em
  minus 0.4em\relax Cambridge University Press, 2008.

\bibitem{neely10}
M.~Neely, \emph{Stochastic network optimization with application to
  communication and queueing systems}.\hskip 1em plus 0.5em minus 0.4em\relax
  Morgan Claypool, 2010.

\bibitem{mckeown99}
N.~McKeown, A.~Mekkittikul, V.~Anantharam, and J.~Walrand, ``Achieving 100%
  throughput in an input-queued switch,'' \emph{{IEEE} Trans. Commun.},
  vol.~47, no.~8, pp. 1260--1267, Aug. 1999.

\bibitem{li11}
Q.~Li and R.~Negi, ``Scheduling in wireless networks under uncertainties: A
  greedy primal-dual approach,'' in \emph{Proc. IEEE Int. Conf. Commun.}, Jun.
  2011, pp. 1--5.

\bibitem{li09}
------, ``Back-pressure routing and optimal scheduling in wireless broadcast
  networks,'' in \emph{Proc. IEEE Global Commun. Conf.}, Dec. 2009, pp. 1--6.

\bibitem{saber11}
A.~Saber and G.~Venayagamoorthy, ``Plug-in vehicles and renewable energy
  sources for cost and emission reductions,'' \emph{{IEEE} Trans. Ind.
  Electron.}, vol.~58, no.~4, pp. 1229--1238, Apr. 2011.

\bibitem{sundstroem10}
O.~Sundstroem and C.~Binding, ``Optimization methods to plan the charging of
  electric vehicle fleets,'' in \emph{Proc. Int. Conf. on Control, Commun. and
  Power Eng.}, 2010.

\bibitem{deilami11}
S.~Deilami, A.~Masoum, P.~Moses, and M.~Masoum, ``Real-time coordination of
  plug-in electric vehicle charging in smart grids to minimize power losses and
  improve voltage profile,'' \emph{{IEEE} Trans. Smart Grid}, vol.~2, no.~3,
  pp. 456--467, Sep. 2011.

\bibitem{turitsyn10}
K.~Turitsyn, N.~Sinitsyn, S.~Backhaus, and M.~Chertkov, ``Robust
  broadcast-communication control of electric vehicle charging,'' in
  \emph{Proc. IEEE SmartGridComm}, Oct. 2010, pp. 203--207.

\bibitem{lasseter11}
R.~Lasseter, ``Smart distribution: Coupled microgrids,'' \emph{Proc. {IEEE}},
  vol.~99, no.~6, pp. 1074--1082, Jun. 2011.

\bibitem{han10}
S.~Han, S.~Han, and K.~Sezaki, ``Development of an optimal vehicle-to-grid
  aggregator for frequency regulation,'' \emph{{IEEE} Trans. Smart Grid},
  vol.~1, no.~1, pp. 65--72, Jun. 2010.

\bibitem{kersting01}
W.~Kersting, ``Radial distribution test feeders,'' in \emph{Proc. IEEE Power
  Eng. Soc. Win. Meet.}, vol.~2, 2001, pp. 908--912.

\bibitem{SCEWebsite}
\BIBentryALTinterwordspacing
S.~C. Edision. (2011) Regulatory information - sce load profiles. [Online].
  Available: \url{http://www.sce.com/AboutSCE/Regulatory/loadprofiles/}
\BIBentrySTDinterwordspacing

\bibitem{EIALoad}
\BIBentryALTinterwordspacing
U.~E.~I. Administration. (2009) Residential average monthly bill by census
  division and state. [Online]. Available:
  \url{http://www.eia.gov/cneaf/electricity/esr/table5.html}
\BIBentrySTDinterwordspacing

\bibitem{CheyVolt}
\BIBentryALTinterwordspacing
U.~E.~P. Agency. (2011) 2011 electric vehicles, fuel economy. [Online].
  Available: \url{http://www.fueleconomy.gov/feg/evsbs.shtml}
\BIBentrySTDinterwordspacing

\bibitem{NHTS}
\BIBentryALTinterwordspacing
P.~Hu and T.~Reuscher. (2004, Dec.) Summary of travel trends. U.S. Department
  of Transportation and Federal Highway Administration. [Online]. Available:
  \url{http://nhts.ornl.gov/2001/pub/STT.pdf}
\BIBentrySTDinterwordspacing

\end{thebibliography}

\begin{IEEEbiographynophoto}{Qiao Li}
  (S'07) received the B.Engg. degree from the Department of Electronics Information Engineering, Tsinghua University, Beijing China, in 2006. He received the M.S. degree from the Department of Electrical and Computer Engineering, Carnegie Mellon University, Pittsburgh, PA USA, in 2008. He is currently a Ph.D. Candidate in the Department of Electrical and Computer Engineering, Carnegie Mellon University. His research interests include distributed algorithms, smart grids technologies, and wireless networking. 
\end{IEEEbiographynophoto}

\begin{IEEEbiographynophoto}{Tao Cui}
  (S'10) was born in Shaanxi Province, China, in 1985. He received the B.Sc. and M.Sc. degrees from Tsinghua University, Beijing, China, and is currently pursuing the Ph.D. degree at Carnegie Mellon University, Pittsburgh, PA. His main research interests include power system computation, protection, analysis, and control.
\end{IEEEbiographynophoto}

\begin{IEEEbiographynophoto}{Rohit Negi}
  received the B.Tech. degree in electrical engineering from the Indian Institute of Technology, Bombay, in 1995. He received the M.S. and Ph.D. degrees from Stanford University, CA, in 1996 and 2000, respectively, both in electrical engineering. Since 2000, he has been with the Electrical and Computer Engineering Department, Carnegie Mellon University, Pittsburgh, PA, where he is a Professor. His research interests include signal processing, coding for communications systems, information theory, networking, cross-layer optimization, and sensor networks. Dr. Negi received the President of India Gold Medal in 1995. 
\end{IEEEbiographynophoto}

\begin{IEEEbiographynophoto}{Franz Franchetti}
  received the Dipl.-Ing. degree and the PhD degree in technical mathematics from the Vienna University of Technology in 2000 and 2003, respectively. Dr. Franchetti has been with the Vienna University of Technology since 1997. He is currently an Assistant Research Professor with the Dept. of Electrical and Computer Engineering at Carnegie Mellon University. His research interests concentrate on the development of high performance DSP algorithms. 
\end{IEEEbiographynophoto}  

\begin{IEEEbiographynophoto}{Marija D. Ili\'c}
  (M'80-SM'86-F'99) is currently a Professor at Carnegie Mellon University, Pittsburgh, PA, with a joint appointment in the Electrical and Computer Engineering and Engineering and Public Policy Departments. She is also the Honorary Chaired Professor for Control of Future Electricity Network Operations at Delft University of Technology in Delft, The Netherlands. She was an assistant professor at Cornell University, Ithaca, NY, and tenured Associate Professor at the University of Illinois at Urbana-Champaign. She was then a Senior Research Scientist in Department of Electrical Engineering and Computer Science, Massachusetts Institute of Technology, Cambridge, from 1987 to 2002. She has 30 years of experience in teaching and research in the area of electrical power system modeling and control. Her main interest is in the systems aspects of operations, planning, and economics of the electric power industry. She has co-authored several books in her field of interest. Prof. Ilic is an IEEE Distinguished Lecturer. 
\end{IEEEbiographynophoto}

\end{document}